\documentclass[12pt, a4paper]{amsart}
\usepackage{amssymb, amsmath, amsfonts, verbatim}
\newtheorem{lemma}[equation]{Lemma}
\newtheorem{theorem}[equation]{Theorem}

\title{New upper bounds on the spreads of the sporadic simple groups}
\author{Ben Fairbairn}
\address{Department of Economics, Mathematics and Statistics, Birkbeck, University of London, Malet Street, London WC1E 7HX}
\email{bfairbairn@ems.bbk.ac.uk}

\begin{document}
\maketitle

\begin{abstract}
We give improved upper bounds on the exact
spreads of many of the larger sporadic simple groups, in some
cases improving on the best known upper bound by several orders of
magnitude.
\end{abstract}



\section{Introduction}
Recall that a group is said to be 2-generated if it is generated
by just two of its elements. Every finite simple group is
2-generated (see \cite{AschbacherGuralnic}) and many authors have
considered the question of how easily a pair of elements
generating a simple group may be obtained. One quantity measuring
this introduced by Brenner and Wiegold in \cite{BrennerWiegold}
and motivated by earlier work of Binder \cite{Binder} is the
concept of the spread of a group.

Let $G$ be a group. We say that $G$ has \emph{spread r} if for any
set of distinct non-trivial elements $X:=\{x_1,\ldots,
x_r\}\subset G$ there exists an element $y\in G$ with the property
that $\langle x_i,y\rangle=G$ for every $1\leq i\leq r$. We say
that this element $y$ is a \emph{mate} of $X$ and that
$G$ has \emph{exact spread s(G):=r} if $G$ has spread $r$ but not $r+1$.

The concept of spread is also of interest to computational group
theorists since it is useful in the the analysis of the celebrated
product replacement algorithm for producing random elements of
groups \cite{Garion}. The concept is also of interest when
studying the generating graph of a group \cite[Section 4]{Maroti}.

The exact spreads of the finite simple groups have been much
studied
\cite{Binder,BrennerWiegold,BreuerGuralnicKantor,GuralnicShalev}.
In particular, bounding the value of the exact spreads of the
sporadic simple groups has recently been investigated by several different authors
\cite{BradleyHolmes,BradleyMoori,GaniefMoori,Woldar} and it is
these cases that we focus on here. More specifically we prove the
following.

\begin{theorem}\label{main}
The exact spreads of the sporadic simple groups are bounded by the values
given in Tables 1 and 2.
\end{theorem}

Most of the bounds listed in Tables \ref{1} and 2 are not new. The upper
bounds given in Table \ref{1} were obtained by Bradley and Holmes
in \cite{BradleyHolmes} using coverings of a group by sets of proper subgroups and as such are unable to handle the sporadic simple groups with very large coverings, which is essentially the larger groups. Our methods are unable to improve upon these bounds.

What is new here are the upper bounds listed in Table 2 for all the groups that the methods of \cite{BradleyHolmes} could not deal with.

The bounds
given in Tables 1 and 2 are, as far as the author is aware, the best
known, accepting that Bradley and Holmes claim that for the groups
they considered ``better results were obtained for some of the
groups in trial runs, but our table gives only the results that
were given by known seeds" \cite[p.138]{BradleyHolmes}.

In Section 2 we will introduce some preliminary ideas that we will
use to prove our bounds in Section 3 in every case aside from the Baby Monster and the
Monster group that we shall deal with separately in Section 4.

\begin{table}\label{1}
\begin{center}
\caption{Bounds on $s(G)$ for the smaller sporadic simple groups. The upper bounds are proved in \cite{BradleyHolmes}. The lower bound for M$_{11}$ is proved in \cite{BradleyHolmes,Breuer,Woldar}. All other lower bounds are proved in \cite{BreuerGuralnicKantor}. (Note that aside from M$_{11}$, M$_{12}$ and J$_2$ the lower bounds stated in \cite[Table 1]{BradleyHolmes} are incorrect.)}
\bigskip

\begin{tabular}{|c|l|c|l|}
\hline $G$&Upper bound&$G$&Upper bound\\
&Lower bound&&Lower bound\\
\hline

M$_{11}$&3&J$_3$&597\\
&3&&76\\
M$_{12}$&9&M$_{24}$&56\\
&3&&11\\
J$_1$&179&McL&308\\
&76&&70\\
M$_{22}$&26&He&1223\\
&20&&198\\
J$_2$&24&Suz&956\\
&5&&40\\
M$_{23}$&8064&Co$_3$&1839\\
&8063&&98\\
HS&33&Fi$_{22}$&186\\
&18&&13\\
\hline
\end{tabular}
\bigskip\small

\end{center}
\end{table}

\begin{table}\label{2}
\caption{The best previous upper bounds for the larger sporadic
groups proved in \cite{BradleyMoori}; the new upper bounds proved here and the
lower bounds proved in \cite{BreuerGuralnicKantor}. Note that the seemingly better lower bound given in
\cite{GaniefMoori} for HN (10999) is incorrect - see \cite[Section
4.7]{BreuerGuralnicKantor}.}
\bigskip

\begin{center}
\begin{tabular}{|c|l|c|l|}
\hline
$G$&Old upper bound&$G$&Old upper bound\\
&New upper bound&&New upper bound\\
&Lower bound&&Lower bound\\
\hline Ru&12990752&Th&103613642531\\
&1252799&&976841774\\
&2880&&133997\\
O'N&5960127&Fi$_{23}$&8853365473\\
&2857238&&31670\\
&3072&&911\\
Co$_2$&5240865&Co$_1$&58021747714\\
&1024649&&46621574\\
&270&&3671\\
HN&229665984&J$_4$&251012689269463297\\
&74064374&&47766599363\\
&8593&&1647124116\\
Ly&112845651178977&Fi$'_{24}$&309163967798745777216\\
&1296826874&&7819305288794\\
&35049375&&269631216855\\
\hline
$\mathbb{B}$&\multicolumn{3}{|l|}{3843675651630431666542962843030}\\
&\multicolumn{3}{|l|}{3843461129719173164840195954999}\\
&\multicolumn{3}{|l|}{174702778623598780219391999999}\\
\hline
$\mathbb{M}$&\multicolumn{3}{|l|}{14587804270839626161268024115186834207944682668030}\\
&\multicolumn{3}{|l|}{5791748068511982636944259374}\\
&\multicolumn{3}{|l|}{3385007637938037777290624}\\
\hline
\end{tabular}
\bigskip\small

\end{center}
\end{table}

\section{Preliminaries}

In this section we shall describe some concepts that will be
useful in proving Theorem \ref{main}.

\subsection{Support}

Let $G$ be a group and let $x\in G^{\#}$ where
$G^{\#}:=G\setminus\{e\}$. We define the \emph{support} of $x$ to be the set
$$supp(x):=\bigcup_{H<G, x\in H}H^{\#}.$$
In other words $y\in supp(x)$ means that $y$ is an element of
$G^{\#}$ that lies in a proper subgroup containing $x$ and so $y$
cannot be a mate for the set $\{x\}$. We extend this to subsets
$X\subset G^{\#}$ as follows:
$$supp(X):=\bigcup_{x\in X}supp(x).$$
If $Y\subset supp(X)$ we say that $X$ \emph{supports} $Y$. In
particular, elements of $supp(X)$ cannot be a mate to $X$.

\subsection{Support classes and characters}\label{class}

Continuing the earlier nomenclature, we say a conjugacy class
$\mathcal{C}\subset G$ is a \emph{support class} if the set
$\mathcal{C}$ supports the set $G^{\#}$. For each of the bounds
that we improve upon here, our improved bound is obtained by
showing that some small conjugacy class of $G$ is a support class.

Note that not every simple group has a support class e.g. the Baby
Monster has no support class as the only maximal subgroups with
elements of order 47 are copies of the Frobenius group 47:23, but
no proper subgroup containing elements of order 23 contains
elements of order 31 - see the list of maximal subgroups given in \cite{baby} (the list given in \cite{ATLAS} is incomplete). Whilst we cannot improve upon the best known bound in this case precisely \emph{the same} methods as the other cases, our approach is not entirely redundant here - see Section \ref{M}.

To obtain a set of elements that has no mate, and thus provide an upper bound on the spread, it is sufficient to take one generating element from each cyclic subgroup generated by an element of a support class. We thus have the following easy lemma.

\begin{lemma}\label{WolLem}
Let $\mathcal{C} \subset G^{\#}$ be a support class, $d:=|\langle g\rangle\cap \mathcal{C}|$ for $g\in\mathcal{C}$. Then
$s(G)+1\leq|\mathcal{C}|/d$.\end{lemma}

Given a conjugacy class $\mathcal{C} \subset G^{\#}$ we define its
\emph{support character} $\chi_{\mathcal{C}}$ to be the sum of the
primitive permutation characters of $G$ that are non-zero on
$\mathcal{C}$. Since the transitive permutation character $1\uparrow_H^G$ is nonzero at a class $\mathcal{C}$ if and only if $\mathcal{C}\cap H\not=\emptyset$ we have the following lemma.

\begin{lemma}\label{char}
A conjugacy class $\mathcal{C} \subset G^{\#}$ is a support class if and
only if $\chi_{\mathcal{C}}(g)>0$ for every $g\in G$.
\end{lemma}

\section{Computing the bounds}

\begin{table}\label{supp}
\caption{Support classes for the sporadic simple groups.}
\bigskip

\begin{center}
\begin{tabular}{|c|c|c|c|c|c|c|c|}
\hline
$G$&Class&$G$&Class&$G$&Class&$G$&Class\\
\hline
M$_{11}$&2A&J$_3$&2A&O'N&2A&Th&2A\\
M$_{12}$&2A&M$_{24}$&2A&Co$_3$&2A&Fi$_{23}$&2A\\
J$_1$&2A&McL&2A&Co$_2$&2B&Co$_1$&2A\\
M$_{22}$&2A&He&2A&Fi$_{22}$&2A&J$_4$&2B\\
J$_2$&2A&Ru&2B&HN&2B&Fi$'_{24}$&2B\\
HS&2A&Suz&3A&Ly&2A&$\mathbb{M}$&2B\\
\hline
\end{tabular}
\end{center}
\end{table}

Our new upper bounds are obtained by finding a small support class using Lemma \ref{char} and then using this to obtain a bound using Lemma \ref{WolLem}.
This is easily done using the GAP algebra system \cite{GAP4}; first by obtaining the primitive permutation characters using standard GAP functions (primarily the GAP character table library and the tables of marks), then obtaining the support characters (if any) using this data and finally by reading off the support class that gives the best bound from the list just obtained. 

The support classes giving the best bound found in this way are given in Table \ref{supp}. For
completeness we give these support classes for each of the
sporadic groups that possesses one, that is every sporadic simple group
aside from the Mathieu group M$_{23}$ and the Baby Monster $\mathbb{B}$.

For example, in several cases (M$_{11}$, J$_1$, M$_{22}$, M$_{23}$, J$_3$,
McL, O'N and Ly) there is only one class of involutions and every
maximal subgroup has even order (see \cite{ATLAS}) so the support character $\chi_{2A}$ is the sum of every primitive permutation character and is
therefore positive on every class. The class 2A is therefore a
support class in these cases by Lemma \ref{char}. Note that whilst the Thompson group, Th, also has only
one class of involutions it also has a class of maximal subgroups
of odd order isomorphic to the Frobenius group 31:15. The
elements of this subgroup can easily be seen to also belong to
other maximal subgroups with structure $2^5.$L$_5(2)$ (see \cite[p.70]{ATLAS}). The support character $\chi_{2A}$ is thus equal to the
sum of every primitive permutation character, aside from the one
defined by the maximal copies of 31:15, and $\chi_{2A}(g)>0$ for every $g\in$ Th. Thus 2A
is a support class by Lemma \ref{char}.

\section{The Baby Monster and the Monster}\label{M}
\subsection{The Baby Monster}

As noted in Section \ref{class} the Baby Monster group $\mathbb{B}$ does not have a support class. However, we can use a union of conjugacy classes instead.

\begin{lemma}
If $\chi:=\chi_{47A}+\chi_{2A}$ then $\chi(g)>0$ for every $g\in\mathbb{B}$.
\end{lemma}

\begin{proof}
Let $g\in\mathbb{B}$.  The group $\mathbb{B}$ has only one class of cyclic subgroups of order 47 so if $o(g)=47$ then we have $\chi(g)>0$. Suppose $o(g)\not=47$. 

Structure constant calculations show that every involution centralizer contains elements of class 2A and so if $g$ is a power of an element of even order we must have $\chi(g)>0$.

The only elements that have yet to be dealt with have order 31. Any such element is contained in a maximal subgroup with structure $[2^{30}]$L$_5(2)$ which can also be shown to contain element of class 2A.
\end{proof}

Combining the above with the natural generalisation of Lemma 2 provides the upper bound given in Table 2. Note that we cannot prove the above in the same computational manner as the results of the previous sections since GAP does not contain all the primitive permutation characters of $\mathbb{B}$ in its libraries.

We further note that only one class of maximal subgroups contains elements of order 47 - copies of the Frobenious group 47:23. To obtain a result like the above we must therefore use either class 23A or class 47A. Using 23A proves an upper bound that is worse than the best previously known upper bound. The above result is therefore the best possible.

\subsection{The Monster}
The Monster group $\mathbb{M}$ requires special attention - the methods of previous sections cannot be applied as easily in this case since, at the time of writing, the maximal subgroups of $\mathbb{M}$, and thus the primitive permutation characters of $\mathbb{M}$, have yet to be classified. 

All not lost! We can still find a support class in this case using information about its conjugacy classes and the known maximal subgroups.

\begin{lemma}\label{lemM}
Class 2B of $\mathbb{M}$ is a support class.
\end{lemma}

\begin{proof}
We aim to show that the primitive permutation characters defined by the known maximal subgroups of $\mathbb{M}$ are sufficient to give us $\chi_{2B}(g)>0$ for every $g\in\mathbb{M}$.

First note that $\mathbb{M}$ has only two classes of involutions (see the character table given in \cite[p.220]{ATLAS}) and so the centralizer of any involution contains 2B elements. Both the 2A and 2B centralizers are known to be maximal \cite[p.228]{ATLAS} and so the sum defining $\chi_{2B}$ must contain both of the permutation characters corresponding to these subgroups. 

Now, let $g\in\mathbb{M}$ and suppose there exist elements $h,k\in\mathbb{M}$ $k\not=1$ such that $g^a=h^b=k$ and $h^c$ is in class 2B for some $a,b,c\in\mathbb{Z}^+$. Then $g,h\in C_{\mathbb{M}}(k)$, which must contain a 2B element. It follows that the sum defining $\chi_{2B}$ must contain the  permutation characters defined by any maximal subgroups containing $C_{\mathbb{M}}(k)$. (For instance, if $g$ is in class 119A then we can find an $h$ in class 14B such that $k:=h^2=g^{17}$ which is in class 7A, so $a=17$, $b=2$ and $c=7$ in this case. A maximal subgroup containing a 7A centralizer will therefore contain elements of class 119A and 2B and so the sum defining $\chi_{2B}$ must contain the permutation character corresponding to such a subgroup.)

Finally, from the fusion maps in the character table we see that the only classes not yet accounted for are the elements of orders 41, 59 and 71. It is known that $\mathbb{M}$ contains maximal copies of 41:40, L$_2$(59) and L$_2$(71) (see for instance \cite[Table 1]{BrayWilson}). Furthermore, it is well known that the product of any two 2A elements of $\mathbb{M}$ has order at most 6. Since each of these subgroups only contain one class of involutions and in each of these subgroups there is a pair of involutions whose product is greater than 6, they must each contain 2B elements. It follows that the sum defining $\chi_{2B}$ must contain the  permutation characters defined by each of these classes of maximal subgroups. 

We thus have that $\chi_{2B}(g)>0$ for all $g\in\mathbb{M}$, so 2B is a support class by Lemma \ref{char}.
\end{proof}

Note we cannot replace 2B by 2A and improve this bound as this would give an `upper bound'
less than the lower bound of \cite{BradleyMoori}.

\section{Concluding remarks}

\begin{enumerate}
\item Our first remark is clear: it would be of great interest to obtain better bounds on the spreads of simple groups, if not determine them all precisely.

In particular, as noted in the MathsciNet review of \cite{BradleyHolmes}, the bounds $3\leq s(\mbox{M}_{12})\leq9$ are tantalizingly close and since M$_{12}$ is such a relatively small and low degree permutation group it seems likely that this particular case is unusually within reach.

\item As Table 3 shows, support classes behave very erratically posing several questions regarding their nature. Which groups possess support classes? Is class 2A a support class infinitely often and if so, for which groups is it a support class? Conversely, among the groups for which 2A is not a support class which other classes are support classes?  Are there groups whose smallest support class has order greater then 3 (class 3A is the smallest sport class for the sporadic Suzuki group) and if so, how large can the order of such a class get? Are there infinitely many groups with a support class of elements of order 3? Of order 4? Of order 5? etc.

It would be of great interest to see answers to all of these questions.

\item There is a more restricted notion of \emph{uniform spread}, where we require the mates of our sets to lie in a single conjugacy class of $G$ independent of the choice of the elements of the sets. In general, spread and uniform spread need not coincide: the group SL$_3$(2) has uniform spread exactly 3 but exact spread 4.

It would be of great interest to determine, or at least bound, the uniform spreads of the sporadic simple groups, which has received much less interest in the literature than exact spreads have \cite{BreuerGuralnicKantor,BurnessGuest}. Clearly the uniform spread is at most the exact spread and so upper bounds, like those proved here, also provide upper bounds on the uniform spread. Otherwise, the only known bounds on the uniform spreads of the sporadic simple groups, as far as the author is aware, are as follows: the uniform spreads of M$_{11}$ and M$_{12}$ are both 3 \cite[Section 5.9 and 5.10]{Breuer}.
\end{enumerate}

\section{Acknowledgments}

The author wishes to express his deepest gratitude to Professor
Jamshid Moori for introducing him to the concept of spread in the
first place and to the referee
for making suggestions that resulted in substantial improvements
to this paper. I am also grateful to Professor Robert Wilson for helpful correspondence regarding the cases of the Baby Monster and Monster groups.


\begin{thebibliography}{14}
\bibitem{AschbacherGuralnic} M. Aschbacher and R.M. Guralnick ``Some Applications of the First
Cohomology Group", J. Algebra, 1984, 90, 446-460

\bibitem{Binder} G.J. Binder ``The two-element bases of the symmetric group", Izv. Vys\v{s}. U\v{c}ebn. Zaved.
Mathematika, 90:9-11, (1970)

\bibitem{BradleyHolmes} J.D. Bradley and P.E. Holmes ``Improved bounds
for the spread of sporadic groups", LMS J. Comput. Math., 10, (2007),
132-140

\bibitem{BradleyMoori} J.D. Bradley and J. Moori ``On the exact spread of sporadic simple
groups", Communications in Algebra, 35(8), 2588-2599, (2007)

\bibitem{BrayWilson} J.N. Bray and R.A. Wilson ``Explicit construction of maximal subgroups of the Monster", J. Algebra, 300, (2006), 835-857

\bibitem{BrennerWiegold} J.L. Brenner and J. Wiegold ``Two-Generator Groups, I", Michigan Math.
J., (1975), 22, 53-64

\bibitem{Breuer} T. Breuer ``GAP Computations Concerning Probabilistic Generation of Finite
Simple Groups", preprint 2007, \texttt{arXiv:0710.3267}

\bibitem{BreuerGuralnicKantor} T. Breuer, R.M. Guralnick and W.M. Kantor
``Probabilistic generation of finite simple groups II", J. Algebra,
320, (2008), 443-494

\bibitem{BurnessGuest} T. Burness and S. Guest ``On the uniform spread of almost simple groups", in preparation

\bibitem{ATLAS} J.H. Conway, R.T. Curtis, S.P. Norton, R.A. Parker and
R.A. Wilson ``An ATLAS of finite groups", (Oxford University
Press, 1985)

\bibitem{GaniefMoori} S. Ganief and J. Moori ``On the spread of the sporadic simple
groups", Communications in Algebra, 29(8), 3239-3255 (2001)

\bibitem{GAP4}
  The GAP~Group, \emph{GAP -- Groups, Algorithms, and Programming,
  Version 4.4.12};
  2008,
  \verb+(http://www.gap-system.org)+

\bibitem{Garion} S. Garion ``Connectivity of the product replacement graph of PSL(2,$q$)", J. Group Theory, 11, (2008), no. 6, 765-777

\bibitem{GuralnicShalev} R.M. Guralnick and A. Shalev ``On the spread of finite simple
groups", Combinatorica, 23, (1), (2003), 73-87


\bibitem{Maroti}  A. Lucchini and A. Mar\'{o}ti ``Some results and questions related to the generating graph of a finite group", Ischia group theory 2008 Proceedings of the conference in group theory, (2009), 183-208

\bibitem{baby} R.A. Wilson ``The Maximal Subgroups of the Baby Monster, I", J. Algebra, 211,
(1999), 1--14

\bibitem{Woldar} A. Woldar ``The exact spread of the Mathieu group
M$_{11}$", J. Group Theory, 10, (2007), 167-171
\end{thebibliography}
\end{document}